\documentclass[a4paper,12pt]{article}

\usepackage{amsthm}
\usepackage{mathrsfs}
\usepackage{amsfonts,amsmath,oldgerm,amssymb,amscd,tikz-cd}
\usepackage{epsfig}
\usepackage{euscript}
\usepackage{amscd}
\usepackage{amsthm}
\usepackage{theoremref}
\usepackage{empheq}

\DeclareMathAlphabet{\mathpzc}{OT1}{pzc}{m}{it}

\newtheorem{thm}{Theorem}[section]

\newtheorem{lem}[thm]{Lemma}
 
\newtheorem{cor}[thm]{Corollary}
\newtheorem{rem}[thm]{Remark}

\newtheorem{question}[thm]{Question}

\newtheorem{definition}{Definition}[section]

\newcommand{\td}{\operatorname{tr.deg}}

\newcommand{\ml}{\operatorname{ML}}

\title{ On characterization of Double Danielewski type algebras }

\author{Parnashree Ghosh$^*$ and Dibyendu Mondal$^{**}$\\
	{\small{\it  $^*$ Department of Mathematics, IISER Pune}}\\ 
	{\small{\it Dr. Homi Bhaba Road, Pune-411008, India}}\\
	{\small{\it  $^{**}$Department of Mathematics, Indian Institute of Technology Indore,}}\\ 
	{\small{\it Simrol, Khandwa Road, Indore-453552, India}}\\
	{\small{\it $^*$e-mail : pmaths@math.tifr.res.in, ghoshparnashree@gmail.com}}\\
	{\small{\it $^{**}$e-mail : mdibyendu07@gmail.com, mdibyendu@iiti.ac.in }}
}

\begin{document}
	\date{}
	\maketitle
	
	\abstract{	Let $k$ be a field. In this paper, we consider Double Danielewski type algebras over an affine factorial $k$-domain $R$.  We observe that this family produces a non-cancellative family of algebras over $R$. Further, when $k$ is a field of characteristic zero, we give a   characterization for an affine algebra to be isomorphic to an algebra of Double Danielewski type.  }
	\noindent
	
	\smallskip
	
	\smallskip
	\noindent
	{\small {{\bf Keywords}. Generalized cancellation problem, Double Danielewski type algebras, Locally nilpotent derivations, Exponential maps, Makar-Limanov invariant. }}
	\smallskip
	
	\noindent
	{\small {{\bf 2020 MSC}. Primary: 14R10, 14R20; Secondary: 13A50, 13B25, 13N15.}}

\section{Introduction}\label{intro}
All rings considered in this paper are commutative Noetherian rings with unity. By $k$ we denote a field.
For an integral domain $A$, the notation $A^{[n]}$ stands for the polynomial ring in $n$ variables over $A$ and $\text{Frac} (A)$ denotes the field of fractions of $A$. 
The group of units of $A$ is denoted by $A^*$. If $T/L$ is a finite field extension then the degree of the extension is denoted by $[T:L]$. A polynomial $f(X) \in A[X]$ is said to be almost monic in $X$, if its highest degree coefficient is a unit in $A$.

We begin with the Generalized Cancellation Problem over integral domains.
\begin{question}\cite{eakin}\thlabel{can que}
Let $R$ be an integral domain and $A$, $B$ be two $R$-domains such that $A^{[1]}=_{R}B^{[1]}$. Does this imply that $A\cong_R B$ ?
\end{question}
When $R=k$ and $A$, $B$ are one dimensional affine domains, then the answer to the above question is affirmative (\cite{AEH}); whereas there are counterexample to the above  question for higher dimensional ($\geq 2$) integral domains (cf. \cite{dani}, \cite{cra1}, \cite{dubo}, \cite{ghosh}). The first family of counterexamples to this question is given by the Danielewski surfaces (\cite{dani}). 
In \cite{mas2}, Masuda has extended the family of Danielewski surfaces to produce non-cancellative algebras over factorial affine domains $R$, containing a field of  characteristic zero. The varieties defined by these algebras are known as {\it hypersurfaces of Danielewski type} (\cite{mas1}). 

In this paper we first exhibit another class of algebras over affine factorial domains $R$, containing fields of arbitrary characteristic, which produce counterexamples to \thref{can que} (cf. \thref{can2}).
In particular, we consider the following family of rings.
\begin{equation}\label{0}
	B_{d,e}^R=\dfrac{R[X,Y,Z,T]}{(X^dY-P(X,Z), X^eT-Q(X,Y,Z))},
\end{equation}
where $R$ is a factorial affine $k$-domain, $d,e\in \mathbb{N}$, $P(X,Z)\in R[X,Z]$ is either monic in $Z$ {or, $P(0,Z)\in R[Z]\setminus R$ is irreducible in $R[Z]$}, and $Q(X,Y,Z)\in R[X,Y,Z]$ is {almost} monic in $Y$, as a polynomial in $\mathrm{Frac}(R)[X,Y,Z]$.  

When $R=k$, the varieties defined by \eqref{0} were introduced by Gupta and Sen in \cite{gupta1} which produce example of a non-cancellative family of surfaces. These surfaces are named as {\it Double Danielewski surfaces}.   
We call algebras of the form $B_{d,e}^R$ as \emph{Double Danielewski type algebras}. We observe that the family $\{B_{d,e}^R\}$ produces a family of counterexamples to \thref{can que} over factorial affine domains containing fields of arbitrary characteristic when $P(X,Z)$ is monic in $Z$ and $\mathrm{deg}_Z(P(X,Z))>1$. (see, \thref{isodd}, \thref{stisodd}). 

 
The fourth section of this paper is devoted to the characterization of the Double Danielewski type algebras over factorial affine domains containing field of characteristics zero.
In \cite{gupta2}, Gupta and Sen characterized all locally nilpotent derivations of double Danielewski surfaces over a field of characteristics zero. Considering their result we give the following criteria for a finitely generated $R$-domain $B$, containing $\mathbb{Q}$, to be isomorphic to a Double Danielewski type algebra $B_{d,e}^R$ as in \eqref{0} (cf. \thref{DD} and \thref{char fac dd}).   

\begin{thm}\thlabel{intro thm1}
Let $R$ be factorial affine domain over a field $k$ of characteristic zero, and let $B$ be a factorial affine $R$-domain. Suppose there exists an irreducible non-trivial locally nilpotent derivation $D:B\rightarrow B$ such that the following hold.
\begin{enumerate}
\item[\rm (i)] $\ker D=R[x]$, where $x\in B$ is transcendental over $R$, 
\item[\rm(ii)] there exists an element $z\in B\setminus xB$ such that $D(z)= x^2 $ and
\item[\rm (iii)] $R[\overline{z}]$ is factorial, where $\overline{z}$ denote the image of $z$ in $B/xB$.
\end{enumerate}
Then $B$ is isomorphic to $R[X,Y,Z,T]/(XY-\varphi(Z), XT-\psi(Y,Z))$ as $R$-algebra, where $\varphi(Z)\in R[Z]\setminus R$, $\psi(Y,Z)\in R[Y,Z]$ is almost monic in $Y$ as a polynomial over $\text{Frac}(R)[Y,Z]$.
\end{thm}

\begin{thm}\thlabel{intro thm2}
Let $R$ be factorial affine domain over a field $k$ of characteristic zero, $K:=\text{Frac}(R)$ and $B$ be a factorial affine $R$-domain. Suppose there exists an irreducible non-trivial locally nilpotent derivation $D$ on $B$ such that the following hold. 
\begin{enumerate}
\item[\rm (i)] $\ker D=R[x]$, where $x\in B$ is transcendental over $R$, 
\item[\rm (ii)] there exists an element $z\in B\setminus xB$ such that $Dz= x^m $ for some $m\geq 2$, and
\item[\rm (iii)] there exists $y\in B\setminus xB$ such that $Dy=x^df(x,z)$ and $B/xB=R[\overline{y},\overline{z}]^{[1]}$, where $0< d < m$, $f(x,z)\in R[x,z]\setminus xR[x,z]$ is  monic in $z$, and $\overline{z},\overline{y}$ denote the images of $z,y$ in $B/xB$.
\item [\rm (iv)] $R[\overline{z}]$ is factorial and $\deg_z f(0,z)=[\text{Frac}(R[\bar{z}]): K]-1$.
\end{enumerate}
Then \[B \cong_{R} R[X,Y,Z,T]/(X^{m-d}Y-F(X,Z), X^dT-G(X,Y,Z))\] where $F(0,Z)\in R[Z]\setminus R$ with $\partial_Z F(X,Z)=f(X,Z)$ and $G(0,Y,Z)$ is almost monic in $Y$ as a polynomial in $\text{Frac}(R)[Y,Z]$.
\end{thm}

These results could be seen as a converse to \cite[Theorem 2.1]{gupta2}. In the next section we mention some basic facts on exponential maps and locally nilpotent derivations on integral domains. 

\section{Preliminaries}\label{prelim}
We begin this section with the definition and some useful properties of exponential maps. 

\begin{definition}
Let $B$ be a $k$-algebra and let $\delta :B\rightarrow B^{[1]}$ be a $k$-algebra homomorphism. For an indeterminate $U$ over $B$, let $\delta_U :B\rightarrow B[U]$ denote the $k$-algebra homomorphism, $\delta_U$ is said to be an \emph{exponential map} if $\delta_U $ satisfies the following properties:
\begin{enumerate}
\item[\rm (i)] For the evaluation map $\epsilon_0: B[U]\rightarrow B$ at $U=0$, the composition $\epsilon_0\delta_U$ is the identity map on $B$. 
\item[\rm (ii)] $\delta_V \delta_U = \delta_{U+V}$, where $\delta_V :B\rightarrow B[V]$ is extended to a homomorphism $\delta_V :B[U]\rightarrow B[V,U]$ by setting $\delta_V(U)=U$.
\end{enumerate}
\end{definition}

The ring of invariant of $\delta$ is denoted by $B^{\delta} :=\{b\in B\, |\, \delta(b)=b \}$. An exponential map $\delta$ is said to be non-trivial if  $B\neq B^{\delta}$. For each $b\in B$, there exists $n_b\in \mathbb{Z}_{\geq 0}$ such that \[ \delta(b) =\sum_{i=0}^{n_b} \delta^{(i)}(b)\, U^i. \] Note that $\{\delta^{(i)}\}_{i\geq 0}$ is a sequence of linear maps on $B$. This sequence is a locally finite iterative higher derivation, which is abbreviated as lfihd. By the \emph{Leibniz rule}  $\delta^{(1)}$ is a $k$-derivation.

\begin{rem}
	\rm Let $B$ be a $k$-domain over a field $k$ of characteristic zero. Let $D$ be a non-trivial locally nilpotent derivation on $B$, then $D$ will induce a non-trivial exponential map $\exp(D)$ on $B$ such that 
	$$
	\exp(D)(b) := \sum_{i=0}^{n_b} \frac{D^ib}{i!} U^i,
	$$ 
	where $n_b$ is the smallest integer such that $D^{n_b+1} (b)=0$.
	Also note that $$\ker D:=\{b \in B \mid Db=0\}=B^{\exp(D)}.$$
\end{rem}

Let $R$ be a $k$-algebra subring of $B$. By EXP$_R(B)$ we denote the set of all exponential maps $\delta$ on $B$ with $R\subseteq B^{\delta}$. Let us denote \[\text{ML}_R (B):=\bigcap_{\delta \in \text{EXP}_R(B)} B^{\delta}\] Then ML$_R(B)$ is a $k$-subalgebra of $B$. Let $\delta: B\rightarrow B[U]$ be an exponential map, then $\delta$ defines a degree function on $B$ denoted by $\deg_{\delta}$. Let $b\in B\setminus \{0\}$ then $\deg_{\delta}(b):=\deg_U(\delta(b))$, and $\deg_{\delta}(0)=-\infty$. Suppose $\delta$ is non-trivial then $a\in B$ is called a \emph{preslice} of $\delta$ if $\deg_{\delta}(a)=$ min$\{\deg_{\delta}(b)\,|\, b\in B\setminus B^{\delta}\}$.

We now note some basic properties of exponential maps (cf. \cite[p. 1291 \& 1292]{crac}).

\begin{lem}\thlabel{prop exp}
Let $B$ be an affine $k$-domain. Suppose there exists a non-trivial exponential map $\delta:B \rightarrow B[U]$. Then the following holds.
 \begin{enumerate}
 \item[\rm (i)] $B^{\delta}$ is a factorially closed subring of $B$. 
 \item[\rm (ii)] $B^{\delta}$ is algebraically closed in $B$. 
 \item[\rm (iii)] If $x\in B$ is such that $\deg_U(\delta(x))$ is of minimal positive degree, and $b$ is the leading coefficient of $U$ in $\delta(x)$, then $b\in B^{\delta}$ and $B[b^{-1}]=B^{\delta}[b^{-1}][x]$. 
 \item[\rm (iv)] $\td_k(B^{\delta})=$ $\td_k(B)-1$.

 \item[\rm (v)] For any multiplicatively closed subset $S$ of $B^{\delta}\setminus \{0\}$, $\delta$ extends to a non-trivial exponential map $S^{-1}\delta$ on $S^{-1}B$ by setting $(S^{-1} \delta)(b/s):= \delta (b)/s$ for $b\in B,\, s\in S$; and the ring of invariant of $S^{-1}\delta$ is $S^{-1}(B^{\delta})$.
 \end{enumerate}
\end{lem}
\begin{rem}
	\em Let $k$ be a field of characteristic zero and $B$ be a $k$-domain. Let $D$ be a non-zero locally nilpotent derivation on $B$ and $\delta$ be the corresponding exponential map $\exp(D)$ on $B$. We then have $\ker D= B^{\delta}$, and hence it satisfies all properties in \thref{prop exp}.
\end{rem}

We end this section by showing that the rings $B_{d,e}^R$ as in \eqref{0} are integral domains. We first state an easy lemma.

\begin{lem}\cite[Lemma 2.4(2)]{Dutta}\thlabel{b}
	Let $A$ be an integral domain and $a,b \in A \setminus \{0\}$. If $b$ is not a zero divisor in $A/(a)$ then the ring $A[T]/(bT-a)$ is an integral domain.
\end{lem}
 Using the above lemma we now prove the following.
\begin{lem}\thlabel{dom lmm}
Let $R$ be a factorial affine $k$-domain. Then the ring $B_{d,e}^R$, as in (\ref{0}), is an integral domain. 
\end{lem}

\begin{proof}
Consider $S:=R[X,Y,Z]/(X^dY-P(X,Z))$. Since $P(0,Z)\in R[Z]\setminus R$, the image of $P(X,Z)$ in $R[X,Z]/ (X)=R[Z]$ is a non-zero divisor, and hence $X^d$ is a non-zero divisor in $R[X,Z]/(P(X,Z))$.
Thus by \thref{b}, $S$ is an integral domain. Let $x,y,z$ denote the residue class of $X,Y,Z$ in $B_{d,e}^R$. Then $S=R[x,y,z]$ is a subring of $B_{d,e}^R$ and $B_{d,e}^R=S[T]/(x^eT-Q(x,y,z))$. 
Now $\frac{S}{(x)}\cong \left(\frac{R[Z]}{(P(0,Z))}\right)[Y]$. Next, we consider two cases.

Suppose $P(X,Z)$ is monic in $Z$ as a polynomial in $R[X, Z]$. Since $Q(0,Y,Z)$ is almost monic in $Y$ as a polynomial in $\text{Frac}(R)[Y,Z]$, it follows that $Q(x,y,z)$ is a non-zero divisor of $S/(x)$. Thus $x^e$ is a non-zero divisor in $S/(Q(x,y,z))$, and hence, by \thref{b}, the assertion follows.

Next, suppose if $P(0,Z)$ is irreducible in $R[Z]$. Then $S/(x)$ is an integral domain. Note that $Q(x,y,z)$ is non-zero, and thus is a non-zero divisor in $S/(x)$. Hence $x^e$ is a non-zero divisor in $S/(Q(x,y,z))$. The assertion follows by applying \thref{b}.
\end{proof}

\section{Non-cancellation of Double Danielewski type varieties}\label{sec dd}
Throughout this section the following notations are fixed. Let $R$ be a factorial affine $k$-domain with field of fraction $K:=\text{Frac}(R)$. Let $X,Y,Z$ and $T$ denote indeterminates over $R$ and
\begin{equation}\label{dd}
B_{d,e}^R:= \dfrac{R[X,Y,Z,T]}{(X^dY-P(X,Z), X^eT-Q(X,Y,Z))},
\end{equation} 	
where $d,e \geq 1$, $r:= \deg_Z P(X,Z) \geq 1$, $s:=\deg_Y Q(X,Y,Z) \geq 1$, $P(X,Z) \in R[X,Z]$ is monic in $Z$, $Q(X,Y,Z)\in R[X,Y,Z]$ is almost monic in $Y$, as a polynomial in $K[X,Y,Z]$, and $d,e,r,s$ satisfy the following conditions
\begin{align}\label{cond}
	&\text{~either~} r \geq 2,\\
	\nonumber& \text{~or~} r=1, s \geq 2 \text{~and~} e \geq 2. 
\end{align}

Let $x,y,z,t$ be the images of $X,Y,Z,T$, respectively, in $B_{d,e}^R$, and $B_{d,e}^{K}:= B_{d,e}^R \otimes_R K$. Note that $B_{d,e}^R$ is an integral domain by \thref{dom lmm}.

We first prove an easy lemma.

\begin{lem}\thlabel{lmm4}
Let $A$ be an affine  $k$-domain and $B$ be a finitely generated $A$-domain. For any multiplicatively closed subset $S$ of $A$, we have $\ml_{S^{-1}A}(S^{-1}B)=S^{-1}\ml_A(B)$. 
\end{lem}

\begin{proof}
Let $\frac{x}{s} \in \ml_{S^{-1}A}(S^{-1}B)$ where $x \in B$ and $s \in S$. Let $\phi \in {\rm EXP}_A(B)$. Since $A\subseteq B^{\phi}$, $\phi$ induces an exponential map $S^{-1} \phi$ on $S^{-1}B$ such that $S^{-1}\phi|_{S^{-1}A}=id_{S^{-1}A}$. Then $(S^{-1}\phi)(\frac{x}{s})=\frac{\phi(x)}{s}=\frac{x}{s}$. Since $B$ is an integral domain, we have $\phi(x)=x$. The choice of the exponential map $\phi \in {\rm EXP}_A(B)$ is arbitrary. Hence it follows that $x \in \ml_A(B)$. Thus $\ml_{S^{-1}A}(S^{-1}B)\subseteq S^{-1}\ml_A(B)$. 

Conversely, let $\frac{y}{t} \in S^{-1}\ml_A(B)$ where $y \in \ml_A(B)$ and $t \in S$. Let $\psi$ be a non-trivial exponential map of $S^{-1}B$, i.e., $\psi \in {\rm EXP}_{S^{-1}A} S^{-1}B$. Since $B$ is a finitely generated $A$-algebra, $\psi$ induces a non-trivial exponential map $\psi_1$ on $B$ such that $\psi_1 |_A= id_A$ and $S^{-1}(B^{\psi_1})=(S^{-1}B)^{\psi}$. As $y \in \ml_A(B)$, $\frac{y}{t} \in S^{-1}(B^{\psi_1})$ and hence we have $\frac{y}{t} \in (S^{-1}B)^{\psi}$. The choice of the non-trivial exponential map $\psi \in {\rm EXP}_{S^{-1}A} S^{-1}B$ is arbitrary. Hence it follows that $\frac{y}{t} \in \ml_{S^{-1}A} S^{-1}B$. Thus $S^{-1}\ml_A(B) \subseteq \ml_{S^{-1}A}(S^{-1}B)$. 
\end{proof}

In \cite[Theorem 3.9]{gupta1}, Gupta and Sen have shown that  $\ml(B_{d,e}^K)=K[x].$ Using their result we prove the following:

\begin{lem}\thlabel{mldd}
 $\ml_R(B_{d,e}^R)=R[x]$.
 \end{lem}
\begin{proof}
 Note that $\ml(B_{d,e}^K)=K[x]$ (cf. \cite[Theorem 3.9]{gupta1}). Now since $\ml_R(B_{d,e}^R)$ is factorially closed in $B_{d,e}^R$, we have $\ml_R(B_{d,e}^R)=S^{-1}\ml(B_{d,e}^R) \cap B_{d,e}^R$ where $S=R \setminus \{0\}$. Thus by \thref{lmm4}, $\ml_R(B_{d,e}^R) = \ml(B_{d,e}^K) \cap B_{d,e}^R=K[x] \cap B_{d,e}^R$. As $B_{d,e}^R \hookrightarrow R[x,x^{-1},z]$, it follows that $\ml_R(B_{d,e}^R)=R[x]$.   
\end{proof}

\begin{rem}\rm 
\thref{mldd} shows that, the Double Danielewski type algebras $B_{d,e}^R$ with conditions  \eqref{cond}, are not $R$-isomorphic to $R^{[2]}$ as $\ml_R(B_{d,e}^R) \neq R$. In particular, for $R=k$, they are not isomorphic to $k^{[2]}$.  
\end{rem}

We now provide some necessary and sufficient conditions for two double Danielewski-type algebras, over affine factorial domains, to be isomorphic. The proof relies on the corresponding necessary condition for double Danielewski algebras over fields established in \cite[Theorem 3.11]{gupta1}. However, the proof of \cite[Theorem 3.11]{gupta1} contained a gap, which has recently been addressed in \cite[Theorem 2.3]{gupta3}, where the authors obtain the result under suitably modified hypotheses on the defining equations of the algebras. 




\begin{thm}\thlabel{isodd}
Let $R$ be a factorial affine $k$-domain and $K=\text{Frac}(R)$. Let \[ B_1:= \dfrac{R[X,Y,Z,T]}{(X^{d_1}Y-P_1(X,Z), X^{e_1}T-Q_1(X,Y,Z))},\] and \[ B_2:= \dfrac{R[X,Y,Z,T]}{(X^{d_2}Y-P_2(X,Z), X^{e_2}T-Q_2(X,Y,Z))}\] be such that $r_i=\deg_Z P_i(X,Z) >1$ and $s_i= \deg_Y Q_i(X,Y,Z)>1$ where $P_i$ is monic in $Z$ and $Q_i$ is almost monic in $Y$, as a polynomial in $K[X,Y,Z]$ for $i=1,2$. 
Suppose that $\rho: B_1 \rightarrow B_2$ is an $R$-isomorphism. Then 

\begin{itemize}
    \item [(i)] $(d_1,e_1,r_1,s_1)=(d_2,e_2,r_2,s_2)=(d,e,r,s)$.

    \item [(ii)] 
    \begin{enumerate}
        \item [(a)] $P_1\bigl(\lambda X, \mu Z+\delta(X)\bigr)=\alpha_1(X,Z)X^d+\mu^r P_2(X,Z),$\\
 for some $\lambda, \mu \in R^*$, $\delta(X) \in R[X]$ and $\alpha_1(X,Z) \in R[X,Z]=R^{[2]}$ with $\mathrm{deg}_Z \alpha_1 <r$.
        
        \item [(b)] $Q_1\bigl(\lambda X, \lambda^{-d}\mu^r Y+ \lambda^{-d} \alpha_1(X,Z), \mu Z+\delta(X)\bigr)
\\ = X^e g_1(X,Y,Z)+\bigl(X^dY-P_2(X,Z) \bigr) g_2(X,Y,Z)+ Q_2(X,Y,Z) g_3(X,Y,Z)\\
 \text{~for some~} g_1,g_2.g_3 \in R[X,Y,Z]=R^{[3]}.$
    \end{enumerate}
\end{itemize}

Conversely, conditions (i) and (ii) are sufficient for $B_1$ and $B_2$ to be $R$-isomorphic.
\end{thm}

\begin{proof}
Let $x,y,z,t$ and $x^{\prime},y^{\prime},z^{\prime},t^{\prime}$ be the images of $X,Y,Z,T$ in $B_1$ and $B_2$, respectively. 
 
 (i)\, First we show that $r_1=r_2$. 
 Since $\rho$ is an $R$-isomorphism it maps $\mathrm{ML}_R(B_1)$ onto $\mathrm{ML}_R(B_2)$. Therefore, by \thref{mldd}, $\rho(R[x]) = R[x']$ and hence $\rho(x)=\lambda x' - \gamma$ for some $\lambda \in R^\ast$ and $\gamma \in R$. Then $\rho$ induces a $K$-isomorphism $\bar{\rho}$ from \[\bar{B}_1=\bigl(K[Y,Z]/(P_1(0,Z),Q_1(0,Y,Z))\bigr)[T]\] onto
    \[\bar{B}_2=\Bigl(K[Y,Z,T]\big/ \bigl(\alpha^{d_2}Y - P_2(\alpha,Z),\; \alpha^{e_2}T - Q_2(\alpha,Y,Z)\bigr) \Bigr)=K^{[1]},\]
    where $\alpha = \lambda^{-1}\gamma$. If possible, suppose $\gamma \neq 0$. Then $\alpha \neq 0$. But in that case
    $\bar{B}_1 \ncong K^{[1]}$, since $\mathrm{deg}_Z P_1(0,Z)>1$. This is a contradiction. Therefore $\gamma=0$, and that implies 
    \begin{equation}\label{b1}
      \rho (x)=\lambda x^{\prime}. 
    \end{equation}
    
Now, $\rho$ extends to the following isomorphism
    \[
    R[x,x^{-1}] \otimes_{R[x]} B_1 = R[x,x^{-1},z]
    \longrightarrow
    R[x',(x')^{-1}] \otimes_{R[x']} B_2= R[x',(x')^{-1},z'].
   \]
    Therefore we have
    \[  u_1 x^i \rho^{-1}(z') = \alpha_1(x)z + \beta_1(x) \] for some $i \ge 0$, $u_1 \in R^\ast$ and $\alpha_1(x),\beta_1(x)\in R[x]$ with $\alpha_1(x)\neq 0$. 
    
   If $i>0$, then $\alpha_1(0)z+\beta_1(0) \in xB_1 \cap R[z]=(P_1(0,z)) R[z]$. This is a contradiction as $\deg_z P_1(0,z)>1$.  Therefore, $i=0$ and thus
   \begin{equation}\label{a1}
        u_1 \rho^{-1}(z') = \alpha_1(x)z + \beta_1(x)
   \end{equation}
   By similar arguments as above, it follows that
\begin{equation}\label{a2}
     u_2 \rho(z) = \alpha_2(x^{\prime})z' + \beta_2(x^{\prime})
\end{equation} where $u_2 \in R^\ast$ and $\alpha_2(x^{\prime}),\beta_2(x^{\prime})\in R[x^{\prime}]$ with $\alpha_2(x^{\prime})\neq 0$.

Now comparing \eqref{a1} and \eqref{a2} and using the fact that $\rho(x)=\lambda x^{\prime}$, we get that
$\alpha_i \in R^{*}$ for $i=1,2$. Thus
\begin{equation}\label{b2}
    \rho(z)=\mu z^{\prime}+\delta(x^{\prime})
\end{equation}
 for some $\mu \in R^{*}$ and $\delta(x^{\prime}) \in R[x^{\prime}]$.

 Now using \eqref{b1} we obtain that, $\rho$ induces an isomorphism 
 \begin{equation}\label{a3}
 \tilde{\rho}: \tilde{B_1}:=\dfrac{R[Y,Z,T]}{\left(P_1(0,Z),Q_1(0,Y,Z) \right)} \longrightarrow \tilde{B_2}:=\dfrac{R[Y,Z,T]}{\left(P_2(0,Z), Q_2(0,Y,Z) \right)}.
 \end{equation}
Let $\tilde{z}$, $\tilde{z}^{\prime}$ denote the images of $Z$ in $\tilde{B_1}$, $\tilde{B_2}$ respectively. Now using \eqref{b2}, and we obtain that $\tilde{\rho}(\tilde{z})=\mu \tilde{z}^{\prime}+\delta(0)$. Now since, $P_1(0,\mu \tilde{z}^{\prime}+\delta(0))=0$, using the fact that $Q_2(0,Y,Z)$ is almost monic in $Y$, as a polynomial in $K[Y,Z]$, it follows that $P_2(0,Z) \mid P_1(0,\mu Z+\delta (0))$ in $R[Z]$ and hence $r_2 \leqslant r_1$. Using similar arguments for the morphism $\tilde{\rho}^{-1}: \tilde{B_2} \to \tilde{B_1}$, we get that $r_1 \geqslant r_2$ and hence the following hold
\begin{equation}\label{c1}
   \boxed{ r_1=r_2=r \text{~(say)~}.}
\end{equation}
Furthermore,
\begin{equation}\label{c2}
    P_1(0,\mu Z+\delta (0))= \eta P_2(0,Z)
\end{equation}
for some $\eta \in R^{*}$.

Next, we now show that $s_1=s_2$. For this let us consider the isomorphism $\tilde{\rho}$ as in \eqref{a3}. This further induces an isomorphism 
$$
\tilde{\rho}^{\prime} : \tilde{B_1} \otimes_R K \to \tilde{B_2} \otimes_R K
$$
such that $\tilde{\rho}^{\prime} (K[\tilde{z}])=K[\tilde{z}^{\prime}]$. Let $\mathcal{M}$ be a maximal ideal of $K[\tilde{z}]$ such that $\tilde{\rho}^{\prime}(\mathcal{M})=\mathcal{N}$ be a maximal ideal of   
$K[\tilde{z}^{\prime}]$. Let $L_1=K[\tilde{z}]/\mathcal{M}$ and $L_2=K[\tilde{z}^{\prime}]/\mathcal{N}$. Suppose $\overline{L_1}$ and $\overline{L_2}$ be the algebraic closure of $L1$ and $L_2$, respectively. Note that $\tilde{\rho}^{\prime}$ induces an isomorphism $\overline{\rho}: \overline{L_1}[Y,T]/(Q_1(0,Y,\tilde{z})) \to \overline{L_2}[Y,T]/(Q_2(0,Y,\tilde{z}^{\prime}))$. Now comparing the number of connected components of $\mathrm{Spec}\left(\overline{L_1}[Y,T]/(Q_1(0,Y,\tilde{z}))\right)$ and $\mathrm{Spec}\left(\overline{L_2}[Y,T]/(Q_2(0,Y,\tilde{z}^{\prime}))\right)$, and their multiplicities we get that 
$$\boxed{s_1=s_2=s \text{~(say)~}.}$$

We now show that $d_1=d_2$ and $e_1=e_2$. Note that the $R$-isomorphism $\rho :B_1\to B_2$ extends to an $K$-isomorphism $\rho_1 : B_1\otimes_R K \to B_2 \otimes_R K$, where \[B_i\otimes_R K=\dfrac{K[X,Y,Z,T]}{(X^{d_i}Y-P_i(X,Z), X^{e_i}T-Q_i(X,Y,Z))}\] for all $i=1,2$.  
Hence from \cite[Theorem 2.3 (II)]{gupta3} it follows that $$\boxed{(d_1,e_1)=(d_2,e_2)=(d,e)\,\text{(say)}}.$$

(ii) Let $E_1:=R[x,z]$ and $E_2:=R[x^{\prime},z^{\prime}]$. From \eqref{b1} and \eqref{b2} it is clear that $\rho(E_1)=E_2$ and hence
$$\rho(x^dB_1 \cap E_1)= (x^{\prime})^d B_2 \cap E_2.$$ From here we get that
\begin{align}
& \nonumber \rho\bigr((x^d, P_1(x,z))E_1\bigl)=\bigl( (x^{\prime})^d, P_2(x^{\prime},z^{\prime})  \bigr)E_2, \\
&\Rightarrow \rho(P_1(x,z))=P_1(\lambda x^{\prime}, \mu z^{\prime}+\delta(x^{\prime}))=\alpha_1(x^{\prime},z^{\prime})(x^{\prime})^d+\beta_1(x^{\prime},z^{\prime})P_2(x^{\prime},z^{\prime}).
\end{align}
Since $P_2(x^{\prime},z^{\prime})$ is monic in $z^{\prime}$ we can assume $\mathrm{deg}_z^{\prime} \alpha_1<r$. Also, as both $P_1(x,z)$ and $P_2(x^{\prime},z^{\prime})$ are monic in $z$ and $z^{\prime}$ of equal degree, it follows that $\beta_1=\mu ^r$. Thus

\begin{equation}\label{P}
    \boxed{\rho(P_1(x,z))=P_1(\lambda x^{\prime}, \mu z^{\prime}+\delta(x^{\prime}))=\alpha_1(x^{\prime},z^{\prime})(x^{\prime})^d+\mu^r P_2(x^{\prime},z^{\prime}).}
\end{equation}
Since, $R[x^{\prime},z^{\prime}]=R^{[2]}$, (ii)-a follows from \eqref{P}. 

Now using \eqref{b1} we have
\begin{equation}\label{y}
    {\rho(y)=\rho(P_1(x,z)/x^d)=\lambda^{-d}\mu^r y^{\prime}+ \lambda^{-d} \alpha_1(x^{\prime},z^{\prime})}
\end{equation}
Let $S_1:=R[x,y,z]$ and $S_2:=R[x^{\prime},y^{\prime},z^{\prime}]$. By \eqref{b1}, \eqref{b2} and \eqref{y}, it follows that $\rho(S_1)=S_2$. Hence we have 
$$
\rho(x^eB_1 \cap S_1)=(x^{\prime})^eB_2 \cap S_2.
$$
From above we get that
\begin{align}\label{rhoQ}
    & \nonumber \rho\bigl((x^e, Q_1(x,y,z)  \bigr)S_1)=\bigl((x^{\prime})^e, Q_2(x^{\prime},y^{\prime},z^{\prime}) \bigr)S_2\\
    & \Rightarrow \bigl( (x^{\prime})^e, Q_1(\lambda x^{\prime}, \lambda^{-d}\mu^r y^{\prime}+ \lambda^{-d} \alpha_1(x^{\prime},z^{\prime}), \mu z^{\prime}+\delta(x^{\prime}) ) \bigr)S_2 = \bigl((x^{\prime})^e, Q_2(x^{\prime},y^{\prime},z^{\prime}) \bigr)S_2
\end{align}
Thus from \eqref{rhoQ} it follows that

\begin{empheq}[box=\fbox]{gather}
\label{Q}
 \nonumber Q_1(\lambda X, \lambda^{-d}\mu^r Y+ \lambda^{-d} \alpha_1(X,Z), \mu Z+\delta(X))\\
\nonumber = X^e g_1(X,Y,Z)+(X^dY-P_2(X,Z) ) g_2(X,Y,Z)+ Q_2(X,Y,Z) g_3(X,Y,Z).
\end{empheq}

for some $g_1,g_2.g_3 \in R[X,Y,Z]=R^{[3]}$. Therefore (ii)-b holds. 

\medskip
Next, we show that the conditions (i) and (ii) are sufficient for $B_1$ and $B_2$ to be isomorphic. Let $(d_1,e_1)=(d_2,e_2)=(d,e)$. 
Define an $R$-algebra homomorphism $\phi:R[X,Y,Z,T] \to B_2$ as follows:
\begin{align*}
    & \phi(X)=\lambda x^{\prime}\\
    & \phi (Z)= \mu z^{\prime} +\delta(x^{\prime})\\
    & \phi (Y) = \lambda^{-d}\mu^r y^{\prime}+ \lambda^{-d} \alpha_1(x^{\prime},z^{\prime})\\
    & \phi (T) = \lambda^{-e} (g_1(x^{\prime},y^{\prime},z^{\prime})+ t^{\prime} g_3(x^{\prime},y^{\prime},z^{\prime})).
    \end{align*}
Now by (ii) it follows that $\phi$ induces $\overline{\phi}: B_1 \to B_2$, as $\bigl(X^dY-P_1(X,Z), X^eT-Q_(X,Y,Z)\bigr) \subset \mathrm{ker}\, \phi$. Now since both $B_1$ and $B_2$ are affine $k$-domains with equal transcendence degree over $k$, it follows that $\overline{\phi}$ is an isomorphism. Hence the result follows. 
\end{proof}

Next, we claim the $R$-isomorphisms of $(B_{d,e}^R)^{[1]}$ and $(B_{d,e-1}^R)^{[1]}$. The proof of the next theorem closely follows the technical details of the proof of \cite[Theorem 3.15]{gupta1}. For the sake of completeness we are including the proof of the theorem.  

\begin{thm}\thlabel{stisodd}
	Let $R$ be a factorial affine $k$-domain and $B_{d,e}^R$ be the integral domain as in (\ref{dd}). Let \[ P(X,Z) =a_0(X)+ a_1(X)Z +\cdots + Z^m \] and \[ Q(X,Y,Z) = b_0(X,Z)+ b_1(X,Z) Y+ \cdots + b_s Y^s \] where $a_i \in R^{[1]}$ for $0 \leq i \leq m$, $b_i \in R^{[2]}$ for $0 \leq i \leq s-1$ and $b_s \in R \setminus \{0\}$. 
We define the following polynomials \[ P^{\prime}(X,Z) = a_1(X)+ 2a_2(X) Z +\cdots +m Z^{m-1} \] and \[ Q^{\prime}(X,Y,Z) = b_1(X,Z) + 2 b_2(X,Z) Y+ \cdots + sb_s Y^{s-1}.\]
Suppose that $e>1$, $\deg_Z P(X,Z)>1$, $(P(0,Z), P^{\prime}(0,Z))=R[Z]$ and \\$(P(0,Z), Q(0,Y,Z), Q^{\prime}(0,Y,Z))$ $ = R[Y,Z]$. Then 
\[ (B_{d,e}^R)^{[1]} \cong_R (B_{d,e-1}^R)^{[1]}.\]

\end{thm}

\begin{proof}
  Let $B:=B_{d,e}^R$.  Recall that $x,y,z$ and $t$ denote the images of $X,Y,Z$ and $T$ in $B$. Let $\phi:B\rightarrow B[U]$ be an exponential map defined on $B$ by 
    \begin{align*}
 	&\phi \mid_R= id_R\\
 	&\phi(x)= x\\
 	&\phi(z)= z+ x^{d+e}U \\
 	&\phi(y)= \frac{P(x,z+x^{d+e}U)}{x^d}=y+x^eU\alpha(x,z,U)\\
 	&\phi(t)= \frac{Q(x,y+x^eU\alpha(x,z,U),z+x^{d+e}U)}{x^e}=t+U\beta(x,y,z,U).
 \end{align*}
 where $\alpha(x,z,U)\in R[x,z,U]$ and $\beta(x,y,z,U)\in R[x,y,z,U]$. 

 Let $A=B[w]=B^{[1]}$. We extend $\phi$ to $A$ by defining $\phi(w)=w-xU$. Let $f=x^{d+e-1}w+z$. Note that $f\in A^{\phi}$.

 Now, similarly as in the proof of \cite[Theorem 3.15]{gupta1}, using the identity 

 $P(x,f)=P(x,z+x^{d+e-1}w)$
 we get $g := y + x^{e-1}\bigl(P'(0,z)\,w + x\,\theta\bigr)$ for some $\theta \in A$ such that 
 \begin{equation}\label{p}P(x,f)=x^dg.\end{equation}
Since $x,f\in A^{\phi}$, it follows that $g \in A^{\phi}$. Again using the identity
$$Q(x,g,f)=Q(x, y + x^{e-1}\bigl(P'(0,z)\,w + x\,\theta\bigr), z+x^{d+e-1}w),$$
we get $h:= P'(0,z)\,Q'(0,y,z)\,w + xt + x\,\rho$ such that 
\begin{equation}\label{r}
Q(x,g,f)=x^{e-1}h.\end{equation} 
Since $x,f,g\in A^{\phi}$, it follows that $h \in A^{\phi}$. 
Now by the given conditions we get $a(Y, Z), b(Y, Z), c(Y, Z) \in R[Y, Z]$ such that
	\begin{equation}
		Q'(0, Y, Z)P'(0, Z)a(Y, Z) + Q(0, Y, Z)b(Y, Z) + P(0, Z)c(Y, Z) = 1.
	\end{equation}
    Now since $P(0,z), Q(0,y,z) \in xA $, we have 
    \begin{equation}\label{q}
    Q'(0,y,z)P'(0,z)a(y,z) +x \eta =1
    \end{equation}
    for some $\eta \in A$. Now using the expresions of $f,g,h$ and \eqref{q} it can be obtained that $w-ha(g,f) \in xA$, and for $v:=w-ha(g,f)/x \in A$ we get $\phi(v)=v-U$.  Thus $A=A^{\phi}[v]=(A^{\phi})^{[1]}$. We now show $A^{\phi} \cong_R  B_{d,e-1}^R$ and hence the result follows.

Let $E = R[x, f, g, h]$ and consider indeterminates $X, F, G, H$ over $k$ so that $R[X, F, G, H] = R^{[4]}$. We define
	\begin{equation}
		E_1 = \frac{R[X, F, G, H]}{(X^d G - P(X, F), X^{e-1} H - Q(X, G, F))} \cong B_{d,e-1}^R.
	\end{equation} 
We first show that $E \cong E_1$. Clearly there exists a surjective $R$-algebra homomorphism $\Phi : R[X, F, G, H] \to E$ such that $\Phi{|_R}=id_R$, $\Phi(X) = x$, $\Phi(F) = f$, $\Phi(G) = g$, and $\Phi(H) = h$. Using equations (\ref{p}) and (\ref{r}), we see that $\Phi$ induces a surjective $R$-algebra map $\overline{\Phi} : E_1 \to E$.
Now note that $A[1/x]=E[1/x][w]=E[1/x]^{[1]}$ and hence $\td_k\, E=\td_k\, A-1$ and further $E_1$ is an integral domain with $\td_k\, E_1=\td_k\, R+2=\td_k\, A-1= \td_k E$. Therefore, $\overline{\Phi}$ must be an isomorphism.

We now show that $A^\phi = E$. Note that $E \subseteq A^\phi$ and since $A[1/x] = E[1/x][w]=E[1/x]^{[1]}$, it follows that $E[1/x] = A^\phi[1/x]$. Therefore, to show that $E = A^\phi$, it is enough to show that $xE = xA^\phi \cap E$. Since $xA^\phi = xA \cap A^\phi$ by \thref{prop exp}(i), it is enough to show that $xE = xA \cap E$, i.e., to show that the kernel of the map $\beta : E \to A/xA$ is $xE$. 
For $a \in A$, let $\tilde{a}$ denote the image of $a$ in $A/xA$. Now

\begin{align}
A/xA &= \frac{R[Y, Z, T, W]}{(P(0, Z), Q(0, Y, Z))} \\
&= \left(\frac{R[Y, Z]}{(P(0, Z), Q(0, Y, Z))}\right)[W, T] \\
&= R[\tilde{z}, \tilde{y}, \tilde{w}, \tilde{t}]
\end{align}

Clearly $\beta(f) = \tilde{z}$, $\beta(g) = \tilde{y}$, and $\beta(h) = P'(0, \tilde{z})Q'(0, \tilde{y}, \tilde{z})\tilde{w}$. 

Now by equation \eqref{q}, $P'(0, \tilde{z})Q'(0, \tilde{y}, \tilde{z}) \in (A/xA)^*$. Hence, $\beta(E) = R[\tilde{z}, \tilde{y}, \tilde{w}]$ and $A/xA = \beta(E)[\tilde{t}] = \beta(E)^{[1]}$. Thus $\dim \beta(E) = \dim(A/xA) - 1$. 

Now since
\begin{equation}
E/xE \cong E_1/xE_1 = R[F, G, H]/(P(0, F), Q(0, G, F)) = \left(\frac{R[F, G]}{(P(0, F), Q(0, G, F))}\right)^{[1]}
\end{equation}

we have $\dim(E/xE) = \dim(A/xA) - 1 = \dim \beta(E)$. Hence kernel of $\beta$ is $xE$. Thus $A^\phi = E \cong_R B_{d,e-1}^R$. 
\end{proof}


Let us consider the following family of affine domains from \eqref{dd}. 
\begin{gather*}
\Omega_3:= \biggl\{ B_{d,e}^R \mid d, e \geq 1, P \text{~is monic in~}Z \text{~as a polynomial in~} R[X,Z] \text{~with~} \deg_Z P>1,\\ Q \text{~is almost monic in~}Y \text{~as a polynomial in~} K[X,Y,Z] \text{~with~} \deg_Y Q >1,\\ (P(0,Z), P^{\prime}(0,Z))=R[Z], \text{~and~} (P(0,Z), Q(0,Y,Z),Q^{\prime}(0,Y,Z))=R[Y,Z] \biggr\}.
\end{gather*}

\begin{cor}\thlabel{can2}
Let $R$ be an affine factorial $k$-domain. Then the family  $\Omega_3$ consists of infinitely many rings which are pairwise non-isomorphic as $R$-algebras and counterexamples to the Generalized Cancellation Problem over $R$.  
\end{cor}

\begin{proof}
Follows from the \thref{isodd} and \thref{stisodd}.
\end{proof}

\section{Characterization of Double Danielewski type algebras}\label{sec char dd}

Throughout this section $R$ is a factorial affine domain over a field $k$ of characteristic zero, $K:=\text{Frac}(R)$ and $B$ is an affine $R$-domain. In this section we give some criteria, in terms of locally nilpotent derivations on $B$, for $B$ to be $R$-isomorphic to an algebra of type \eqref{0}.
The following theorem of Gupta and Sen classifies all locally nilpotent derivations of $B^k_{d,e}$ under suitable conditions.

\begin{thm}\cite[Theorem 2.1]{gupta2}\thlabel{chara lnd}
	Let $$B^k_{d,e}=\dfrac{k[X,Y,Z,T]}{(X^dY-P(X,Z),X^eT-Q(X,Y,Z))},$$ where $r:=\deg_ZP(X,Z)$ and $s:=\deg_YQ(X,Y,Z)$. Suppose that $P$ and $Q$ are monic in $Z$ and $Y$, respectively, with $r$ and $s$ satisfying one of the conditions in \eqref{cond}. Then any locally nilpotent derivation of $B^k_{d,e}$ is of the form $fD$, where $f\in k[x]$ and $D$ is given by 
	\begin{align*}
			&D(x)= 0, \\
			&D(z)= x^{d+e},\\ 
			&D(y)=  \partial_zP(x,z)\, x^e,  \text{~and~} \\
			&D(t)= \partial_y Q(x,y,z)\,\, \partial_z P(x,z) + \partial_z Q(x,y,z)\, x^d.
	\end{align*}
\end{thm}



We now prove \thref{intro thm1} and \thref{intro thm2}. These results could be seen as converse of \thref{chara lnd}.  
First we state  and prove some preparatory lemmas. 

\begin{lem}\cite[Lemma 3.9]{mas2}\thlabel{lmm1}
Let $R$ be a factorial affine $k$-domain. Let $B$ be an affine $R$-domain generated over $R$ by the elements $x,y,z \in B$. Suppose that $x$ and $z$ are algebraically independent over $R$ and the generators satisfy $x^my-F(x,z)$, where $m\geq 1$ and $F(x,z)\in R[x,z]$ is such that $F(0,z)\in  R[z]\setminus R$. Then $B$ is isomorphic to $R[X,Y,Z]/(X^mY-F(X,Z))$ as $R$-algebra. 
\end{lem}

\begin{lem}\thlabel{lmm3}
Let $R$ be a factorial affine $k$-domain. Let $B$ be an affine $R$-domain generated over $R$ by the elements $x,y,z,t \in B$ with $\td_kB=\td_kR + 2$. Suppose that the generators satisfy relations $x^dy-F(x,z)$ and $x^et=G(x,y,z)$, where $e,d\geq 1$, $F(x,z)\in R[x,z]$ is such that $F(0,z)\in  R[z]\setminus R$ and is irreducible in R[z], and $G(x,y,z)\in R[x,y,z]$ is such that $G(0,y,z)$ is almost monic in $y$ as a polynomial in $K[y,z]$. Then $B$ is isomorphic to $R[X,Y,Z,T]/(X^dY-F(X,Z),X^eT-G(X,Y,Z))$ as $R$-algebra. 
\end{lem}

\begin{proof}
Let $\Phi: R[X,Y,Z,T]\rightarrow B$ be the $R$-algebra homomorphism defined by $\Phi(X)=x$, $\Phi(Y)=y$, $\Phi(Z)=z$ and $\Phi(T)=t$. Then $(X^dY-F(X,Z),X^eT-G(X,Y,Z))\subseteq $ $\ker\Phi$. Let $S:=R[X,Y,Z,T]/(X^dY-F(X,Z),X^eT-G(X,Y,Z))$. Note that $S$ is an integral domain (cf., \thref{dom lmm}). Then the induced homomorphism $\overline{\Phi}: S\rightarrow B$ is a surjective $R$-algebra homomorphism between integral domains with same transcendence degree over $k$. Hence $\overline{\Phi}$ is an $R$-algebra isomorphism.  
\end{proof}

The following lemma is due to Masuda and a proof is contained in \cite[Theorem 4.1]{mas2}. However, for the sake of completeness, we are adding the proof here.
\begin{lem}\thlabel{slice}
Let $k$ be a field of characteristic zero, $R$ be a factorial affine $k$-domain and $B$ be a factorial affine $R$-domain. Suppose there exists a non-trivial locally nilpotent derivation $D$ on $B$ such that the following holds.
\begin{enumerate}
\item[\rm (i)] $\ker D=R[x]$, where $x\in B$ is transcendental over $R$,
\item[\rm (ii)] there exists $z\in B\setminus xB$ such that $Dz=x^m$ for some $m\in \mathbb{Z}_{\geq 1}$.
\end{enumerate}
If $B\neq R[x,z]$, then $R[z]\cap xB=\varphi(z)R[z]$ where $\varphi(z)$ is an irreducible element of $R[z]$ with $\deg_z \varphi(z) \geq 1$.
\end{lem}

\begin{proof}
From \thref{prop exp}(i), it follows that $x$ is an irreducible element of $B$, and hence a prime element of $B$. Also, by \thref{prop exp}(iii), $B\subseteq B[x^{-1}]=R[x,x^{-1},z]$. Thus $x$ and $z$ are algebraically independent over $R$, and $R[z]\cap xB$ is a prime ideal of $R[z]$. Since $R[x,z]\subsetneq B$, there exists an element $b\in B$ such that $x^nb=F(x,z)$ where $F(x,z)\in R[x,z]\setminus xR[x,z]$ and $n\geq 1$. Write $F(x,z) = xG(x,z) + H(z)$, where $G(x,z)\in R[x,z]$ and $H(z)\in R[z]$. Then $H(z)$ is non-zero and belongs to $R[z]\cap xB$, and hence $R[z]\cap xB \neq (0)$. Let us consider the ideal $K[z]\cap x(B\otimes_R K)$ of $K[z]$. Since $K[z]$ is a principal ideal domain and $K[z]\cap x(B\otimes_R K)$ is a non-zero prime ideal of $K[z]$, we have $K[z]\cap x(B\otimes_R K)=\varphi(z)K[z]$ for some non-constant irreducible element $\varphi(z)$ of $K[z]$ with $\deg_z \varphi(z) \geq 1$. As $R$ is factorial, we can assume that $\varphi(z)$ is an irreducible element of $R[z]$. Also, $\varphi(z)K[z]\cap R[z]=\varphi(z)R[z]$. Hence, $R[z]\cap xB=\varphi(z)R[z]$ where $\varphi(z)$ is an irreducible element of $R[z]$ with $\deg_z\varphi(z)\geq 1$.
\end{proof}


We now prove two technical lemmas which will be crucial to prove our main results on characterization of $B_{d,e}^R$.

\begin{lem}\thlabel{lem2}
	Let $k$ be a field of characteristic zero, $R$ be a factorial affine $k$-domain and $B$ be a factorial $R$-domain. Suppose there exists an irreducible locally nilpotent derivation $D$ on $B$ such that the following hold.
	\begin{itemize}
		\item [\rm (i)] $\ker D=R[x]$, where $x$ is transcendental over $R$,
		
		\item [\rm (ii)] there exist $z \in B \setminus xB$ such that $Dz=x^m$ for some integer $m>0$, and 
		
		\item [\rm (iii)] there exists an element $y\in B$ such that $Dy \in xB$ and $z,y$ are algebraically independent over $R$.
	\end{itemize} 
If $R[\overline{z}]$ is factorial then $R[y,z] \cap xB=(\varphi(z), \psi(y,z)) R[y,z]$ where $\overline{z}$ denote the image of $z$ in $B/xB$, $\varphi(z)$ is an irreducible polynomial in $R[z]$ and $\psi(y,z)$ is almost monic in $y$ as a polynomial in $K[y,z]$.
\end{lem}
\begin{proof}
We consider the prime ideal $\mathfrak{p}:= R[y,z]\cap xB$ of $R[y,z]$. By \thref{slice}, $R[z] \cap xB=\varphi(z)R[z]$, for some irreducible polynomial in $R[z]$.
 Note that $\varphi(z)R[z]=R[z]\cap xB \subseteq R[y,z]\cap 
  xB$,  and hence $\mathfrak{p}$ is a non-zero prime ideal of $R[y,z]$. Let $\overline{\mathfrak{p}}$ denote the prime ideal $\mathfrak{p}/\varphi(z)R[y,z]$ of $R[\overline{y}, \overline{z}]$. Now since $D$ is irreducible, it induces a non-trivial locally nilpotent derivation $\overline{D}$ on $B/xB$ such that $R[\overline{y}, \overline{z}] \subseteq \ker \overline{D}$. Since $\td_R (\ker \overline{D})=0$, it follows that $\overline{y}$ is algebraic over $R[\overline{z}]$.
Now if $\mathfrak{p} =(\varphi(z))R[y,z]$, then it follows that $\overline{y}$ is algebraically independent over $R[\overline{z}]$, which is a contradiction. Therefore, $(\varphi(z))R[y,z] \subsetneq \mathfrak{p}$, i.e., $\overline{\mathfrak{p}} \neq (0)$. 

Consider the ideal $\mathfrak{p}_K = K[y,z]\cap x(B \otimes_R K)= \mathfrak{p} \otimes_R K$ of $K[y,z]$. Since $x$ is algebraically independent over $R$, and $\ker D=R[x]$ is factorially closed in $B$ it follows that $\mathfrak{p} \cap R=(0)$ and thus $\mathfrak{p}_K$ is a prime ideal in $K[y,z]$. Now consider the prime ideal $\overline{\mathfrak{p}_K}= \mathfrak{p}_K/ \varphi(z) K[y,z]= \overline{\mathfrak{p}} \otimes_R K$ of $K[\overline{z}][y]$. Since $\mathfrak{p} \neq 0$, it follows that $\mathfrak{p}_K \neq 0$ and hence $\overline{\mathfrak{p}_K} =\psi(y,\overline{z})K[\overline{z}][y]$ where $\psi(y,\overline{z})$ is an irreducible polynomial of $K[\overline{z}][y]$ with $\deg_y \psi(y,\overline{z}) \geq 1$. We can assume that $\psi(y,\overline{z})$ is an element of $R[\overline{z}][y]$ and the highest power of $y$ has coefficient in $R$. Since $R[\overline{z}]$ is factorial it follows that $\psi(y,\overline{z})K[\overline{z}][y]\cap R[\overline{z}][y]=\psi(y,\overline{z})R[\overline{z}][y]$, and hence we deduce that $\overline{\mathfrak{p}}=\psi(y,\overline{z})R[\overline{z}][y]$.
Thus the assertion follows.
\end{proof}

\begin{lem}\thlabel{lem1}
	Let $R$ be an integral domain, $B$ be an $R$-domain and $n\geq 2$ be an integer. Let $a, b_1, \ldots, b_n \in B \setminus R$ be such that the following conditions hold.
	\begin{itemize}
		\item [\rm (i)] $B \subseteq R[a,a^{-1}, b_1,\ldots,b_n]$, 
		
		\item [\rm (ii)] $S:=R[b_1,\ldots,b_{n-1}]=R^{[n-1]}$ with \\$S \cap aB=\left(\varphi_1(b_1), \varphi_2(b_1,b_2),\ldots,\varphi_{n-1}(b_1,\ldots,b_{n-1})\right)S$ such that 
		$$
		\varphi_i(b_1,\ldots,b_i)=a^{l_i}g_i(a,b_1,\ldots,b_{i+1}) \in R^{[i]} \,\, \text{for all}\,\,  1 \leqslant i \leqslant n-1, 
		$$
		where $g_i \in R^{[i+2]}$ and integers $l_1,\ldots,l_{n-1} >0$, and 
		
		\item [\rm (iii)] $\overline{b_n}$ is algebraically independent over $R[\overline{b_1}, \ldots, \overline{b_{n-1}}]$, where $\overline{b_1}, \ldots, \overline{b_n}$ denote the images of $b_1,\ldots,b_n$ in $B/aB$.
	\end{itemize}
Then $B=R[a,b_1,\ldots,b_n]$.
\end{lem}
\begin{proof}
  Let $b \in B$. By (i) there exists $m \geqslant0$, such that \begin{equation}\label{a}
  	a^mb=f(a,b_1,\ldots,b_n) 
  \end{equation} 
  for some $f \in R^{[n+1]}$. By induction on $m$ we show that $b \in R[a,b_1,\ldots,b_n]$. Since $b \in B$ is arbitrary, thus our assertion will follow. 
  
  If $m=0$, then we are already done. Therefore, we assume that $m>0$.
  
  Suppose that if $a^{m-1}b \in R[a,b_1,\ldots,b_n]$, then $b \in R[a,b_1,\ldots,b_n]$. 
  
  Now from \eqref{a} we have
  \begin{equation}\label{c}
  a^mb= a f_1(a,b_1,\ldots,b_n) + \sum_{i=0}^{s} h_i(b_1,\ldots,b_{n-1}) b_n^i.
\end{equation}
From \eqref{c} it follows that $\sum_{i=0}^{s} h_i(\overline{b_1},\ldots,\overline{b_{n-1}}) \overline{b_n}^i=0$ in $B/aB$. Now by (iii) we have $h_i(\overline{b_1},\ldots,\overline{b_{n-1}})=~0$ in $B/aB$ for $0 \leqslant i \leqslant s$, i.e., $h_i(b_1,\ldots,b_{n-1}) \in S \cap aB$. Now since $B$ is an integral domain, by (ii) we have $a^{m-1}b \in R[a,b_1,\ldots,b_n]$. Thus the desired assertion follows by induction hypothesis. 
\end{proof}

We now prove our main results on characterization of the Double Danielewski type algebras $B_{d,e}^R$.
\begin{thm}\thlabel{DD}
Let $k$ be a field of characteristic zero, $R$ be a factorial affine $k$-domain, $K:=\text{Frac}(R)$ and $B$ be a factorial affine $R$-domain. Suppose there exists an irreducible locally nilpotent derivation $D$ on $B$ such that the following holds.
\begin{enumerate}
\item[\rm (i)] $\ker D=R[x]$, where $x\in B$ is transcendental over $R$, 
\item[\rm(ii)] there exists an element $z\in B\setminus xB$ such that $Dz= x^2$ and
\item[\rm (iii)] $R[\overline{z}]$ is factorial, where  $\overline{z}$ is the image of $z$ in $B/xB$.
\end{enumerate}
Then $B$ is isomorphic to $R[X,Y,Z,T]/(XY-\varphi(Z), XT-\psi(Y,Z))$ as $R$-algebra, where $\varphi(Z)\in R[Z]\setminus R$, $\psi(Y,Z)\in R[Y,Z]$ is almost monic in $Y$ as a polynomial over $K[Y,Z]$.
\end{thm}

\begin{proof}
By \thref{prop exp}(i), $x$ is irreducible, and hence prime in $B$. From \thref{prop exp}(iii), we get that $B\subseteq R[x,x^{-1},z]$. Thus $x$ and $z$ are algebraically independent over $R$. Suppose $B=R[x,z]=R^{[2]}$, then $B$ is isomorphic to $R[X,Y,Z,T]/(XY-\varphi(Z),XT-\psi(Y,Z))$ as an $R$-algebra, where $\varphi(Z)=Z\in R[Z]\setminus R$ and $\psi(Y,Z)=Y\in R[Y]\setminus R$. Henceforth, we only consider the case when $R[x,z]\subsetneq B$. We use the following notation for rest of the proof: $\overline{B}:=B/xB$, $\overline{D}$ denotes the induced locally nilpotent derivation on $\overline{B}$ by $D$, and $\overline{b}$ denotes the residue class of $b\in B$ in $\overline{B}$. Note that
\begin{equation}\label{e1}
R[z]\cap xB=\varphi(z)R[z]
\end{equation}
where $\varphi(z)$ is an irreducible polynomial of $R[z]$ with $\deg_z \varphi(z)\geq 1$ (c.f. \thref{slice}). Since $\varphi(z)\in xB\setminus \{0\}$, we can write 
\begin{equation}\label{e2}
\varphi(z)=x^{l_1}y
\end{equation}
where $y\in B\setminus xB$ and $l_1\geq 1$. Applying $D$, it follows that $\partial_z \varphi(z) x^2 = x^{l_1} (Dy)$. Suppose $l_1 > 2$, then $\partial_z \varphi(z)=x^{l_1-2} (Dy)\in xB$. Hence $\partial_z \varphi(z) \in R[z]\cap xB = \varphi(z) R[z]$, which is a contradiction as $\deg_z (\partial_z \varphi(z)) < \deg_z \varphi(z)$. Thus $l_1 \leq 2$, and 
\begin{equation}\label{e3}
Dy = x^{2-l_1}\, (\partial_z \varphi(z)). 
\end{equation}
 
 Suppose $l_1=2$. Then $Dy = \partial_z \varphi(z)$. Since $\partial_z \varphi(\overline{z})\neq 0$ in $\overline{B}$ and the subring $R[\overline{z}]$ of $\overline{B}$ is contained in $\ker \overline{D}$, the element $\overline{y}\notin \ker \overline{D}$ and hence it is algebraically independent over $R[\overline{z}]$. Therefore, $\overline{D}$ is a non-trivial locally nilpotent derivation on $\overline{B}$. Now, since $B \subseteq R[x,x^{-1},z,y]$, $x^2y=\varphi(z)$ and  $\overline{y}$ is algebraically independent over $R[\overline{z}]$, by \thref{lem1}, it follows that $B=R[x,y,z]$. The generators of $B$ satisfy the relation $x^2y=\varphi(z)$ such that $\varphi(z)\in R[z]\setminus R$. From \thref{lmm1}, $B$ is $R$-isomorphic to $R[X,Y,Z]/(X^2Y-\varphi(Z))$, where $\varphi(Z)\in R[Z]\setminus R$.  That means $B \cong_R R[X,Y,Z,T]/(XY-\varphi(Z), XT-\psi(Y,Z))$ where $\psi(Y,Z)=Y$.

Now we assume that $l_1=1$, that is
\begin{equation}\label{1}
	\varphi(z)=xy.
\end{equation}
Then 
\begin{equation}\label{2}
	Dy = x (\partial_z \varphi(z)), 
\end{equation}
and hence $\overline{y}\in \ker \overline{D}$. By \eqref{1}, it follows that $y$ and $z$ are algebraically independent over $R$.  Since conditions (1) and (2) hold and $R[\overline{z}]$ is factorial, by \thref{lem2} we have
 \begin{equation}\label{d}
 	\mathfrak{p}=R[y,z] \cap xB=(\varphi(z), \psi(y,z))R[y,z], 
 \end{equation}
 where $\psi(y,z)\in R[y,z]$ is almost monic in $y$ as a polynomial in $K[y,z]$.

Now let us write 
\begin{equation}\label{e}
	\psi(y,z)=x^{n}t,
\end{equation}
where $t\in B\setminus xB$ and $n\geq 1$. Applying $D$ on \eqref{e}, we get that
\[x^nD(t)=x^2 (\partial_z \psi(y,z)) + x (\partial_y \psi(y,z)) (\partial_z \varphi(z)). \] 

Suppose $n>1$. Then $(\partial_y \psi(y,z)) (\partial_z \varphi(z)) \in R[y,z]\cap xB=(\varphi(z),\psi(y,z))R[y,z]$ and therefore, $\partial_y \psi(y,\bar{z})$ $ \partial_z \varphi(\bar{z}) \in \psi(y,\bar{z})K[\bar{z}][y]$. Since $\deg_y \psi(y,\bar{z}) \geq 1$ and $\partial_y \psi(y,\bar{z})$ is non-zero in $R[\bar{z}][y]=R[\bar{z}]^{[1]}$, we get a contradiction as $\deg_y (\partial_y \psi(y,\bar{z})) < \deg_y \psi(y,\bar{z})$. Thus $n=1$ and
\begin{equation}\label{3}
	 \psi(y,z)=xt 
\end{equation}
where $t\notin xB$. 
Now, from \eqref{3}, it follows that $D(t)=x (\partial_z \psi(y,z)) + (\partial_y \psi(y,z)) (\partial_z \varphi(z))$, which implies that $\overline{t}\notin \ker (\overline{D})$. 

Now, since $R[\overline{z},\overline{y}]\subseteq \ker (\overline{D})$, it follows that $\overline{t}$ is algebraically independent over $R[\overline{y}, \overline{z}]$ (cf. \thref{prop exp}(ii)).
Note that since $B \subseteq R[x,x^{-1},y,z,t]$, $\varphi(z)=xy$, $\psi(y,z)=xt$ and $\overline{t}$ is algebraically independent over $R[\overline{y}, \overline{z}]$, by \thref{lem1}, it follows that $B=R[x,y,z,t]$.   
The generators of $B$ satisfy the equations \eqref{1} and \eqref{3}, the assertion follows from \thref{lmm3}. 
\end{proof}


We now prove \thref{intro thm2} of Section \ref{intro}.

\begin{thm}\thlabel{char fac dd}
Let $R$ be a factorial affine domain over a field $k$ of characteristic zero, $K:= \text{Frac}(R)$ and $B$ be a factorial affine $R$-domain. Suppose there exists an irreducible locally nilpotent derivation $D$ on $B$ such that the following hold.
\begin{enumerate}
\item[\rm (i)] $\ker D=R[x]$, where $x\in B$ is transcendental over $R$, 
\item[\rm (ii)] there exists an element $z\in B\setminus xB$ such that $Dz= x^m $ for some $m\geq 2$, and
\item[\rm (iii)] there exists $y\in B\setminus xB$ such that $Dy=x^df(x,z)$ and $B/xB=R[\overline{y},\overline{z}]^{[1]}$, where $0< d < m$, $f(x,z)\in R[x,z]\setminus xR[x,z]$ is monic in $z$, and $\overline{z},\overline{y}$ denote the images of $z,y$ in $B/xB$.
\item [\rm (iv)] $R[\overline{z}]$ is factorial and $\deg_z f(0,z)=[\text{Frac}(R[\bar{z}]): K]-1$.
\end{enumerate}
Then \[B \cong_{R} R[X,Y,Z,T]/(X^{m-d}Y-F(X,Z), X^dT-G(X,Y,Z))\] where $F(0,Z)\in R[Z]\setminus R$ with $\partial_Z F(X,Z)=f(X,Z)$ and $G(0,Y,Z)$ is almost monic in $Y$ as a polynomial in $K[Y,Z]$.
\end{thm}

\begin{proof}
Note that $x$ is prime in $B$, $B\subseteq R[x,x^{-1},z]$, and $x$ and $z$ are algebraically independent over $R$. As in \thref{DD}, we only consider the case $B\subsetneq R[x,z]$. We use the following notations for the rest of this proof: $\overline{B}:=B/xB$, $\overline{D}$ denotes the induced non-zero locally nilpotent derivation on $\overline{B}$ by $D$, and $\overline{b}$ denotes the residue class of $b\in B$ in $\overline{B}$. 

Let $y \in B$ as in (iii). Since $B \subsetneq R[x,x^{-1},z]$, there exists $n\geq 0$ such that
\begin{equation}\label{f}
x^ny=F(x,z)
\end{equation} 
for some $F(x,z)\in R[x,z]\setminus xR[x,z]$. Applying $D$ on the above equation and using (iii), we have 
\begin{equation}\label{g}
	x^{n+d}f(x,z)=x^n Dy  =x^m (\partial_z F(x,z)).
\end{equation}
If $m> n+d$, then $f(x,z)\in xR[x,z]$, which is a contradiction. Now, suppose $m< n+d$, then $\partial_z F(x,z) \in xR[x,z]$. Let $F(x,z)=\sum_{i=0}^s H_i(x)z^i$ for some $s\in \mathbb{Z}_{\geq 0}$ and $H_i(x)\in R[x]$ for all $0\leq i\leq s$. 
Since $R[x,z]= R^{[2]}$, $H_i(x)\in xR[x]$ for all $1\leq i\leq s$. Now from \eqref{f},  it follows that $H_0(x)\in xB \cap R[x]$. Since $R[x]$ is factorially closed in $B$, by \thref{prop exp}(ii) we have $H_0(x) \in xR[x]$.
But then $F(x,z)\in xR[x,z]$, which is a contradiction. 
Thus $n=m-d$, and hence
\begin{equation}\label{4}
	x^{m-d}y=F(x,z).
\end{equation}
 Now from \eqref{g} it follows that 
 \begin{equation}\label{f1}
 	\partial_z F(x,z)=f(x,z)
 \end{equation} 
 and therefore, monicness of $f$ in $z$ implies that $F$ is almost monic in $z$, i.e., the coefficient of the highest degree term of $F$ in $z$ is unit in $R$. Also by \thref{slice}, 
\begin{equation}\label{e7}
R[z]\cap xB=\varphi(z)R[z]
\end{equation}
where $\varphi(z)$ is an irreducible polynomial of $R[z]$ with $\deg_z \varphi(z)\geq 1$.
 We now show that $F(0,z)=\alpha \varphi(z)$ for some $\alpha \in R^{*}$. Since $R[\overline{z}]=R[z]/\varphi(z)R[z]$, by condition (iv), we have $\deg_z f(0,z)=\deg_z \varphi(z) -1$. Therefore, $\deg_z F(0,z)=\deg_z \varphi(z)$. Also by \eqref{4} we have,  $F(0,z)\in R[z]\cap xB=\varphi(z)R[z]$. Hence 
\begin{equation}\label{e8}
F(0,z)=\beta \varphi(z)
\end{equation}
for some $\beta \in R\setminus \{0\}$. By almost monicness of $F(0,z)$, we get that $\beta \in R^*$. Thus 
\begin{equation}\label{e19}
\varphi(z)=x^{m-d}\beta^{-1} y+ x H_1(x,z)
\end{equation}
for some $H_1(x,z)\in R[x,z]$. 

Since $R[\overline{y},\overline{z}]^{[1]}=\overline{B}$ and $R[\overline{y},\overline{z}]\subseteq \ker \overline{D}$, it follows from \thref{prop exp}(ii) that $\ker \overline{D}=R[\overline{y},\overline{z}]$. 
Also, $y$ and $z$ are algebraically independent over $R$. 
We now consider the prime ideal $\mathfrak{p}:=R[y,z]\cap xB$ of $R[y,z]$. Since conditions (i), (ii) and (iii) are satisfied and $R[\overline{z}]$ is factorial, by \thref{lem2} we have $$\mathfrak{p}=(\varphi(z),\psi(y,z))R[y,z]$$
 where $\psi(y,z)\in R[y,z]$ is almost monic in $y$ as a polynomial in $K[y,z]$. 
Furthermore, $\psi(y,\overline{z})$ is an irreducible element of $R[\overline{z}][y]$ with $\deg_y \psi(y,\overline{z})\geq 1$. Since $\psi(y,z)\in xB$, we have 
\begin{equation}\label{e11}
\psi(y,z)=x^{t_1}v_1
\end{equation}
where $v_1\in B\setminus xB$ and $t_1\geq 1$. Applying $D$ on \eqref{e11} we obtain the following 
\begin{equation}\label{e12}
x^{t_1}Dv_1=x^m (\partial_z \psi(y,z)) + x^d (\partial_y \psi(y,z)) (\partial_z F(x,z)). 
\end{equation}
If $t_1> d$, from \eqref{f1} and \eqref{e8}, it follows that $(\partial_y \psi(y,z)) (\partial_z \varphi(z)) \in xB \cap R[y,z]=(\varphi(z), \psi(y,z))R[y,z]$. But then we have a contradiction as $\deg_z \partial_z \varphi(z)< \deg_z \varphi(z)$ and $\deg_y \partial_y \psi(y,\overline{z}) < \deg_y \psi(y,\overline{z})$, where $\overline{z}$ denote the image of $z$ in $R[y,z]/\varphi(z)R[y,z]$.
 Hence $t_1\leq d$. Suppose $t_1<d$, then $\overline{v_1}\in \ker \overline{D}$.
 Therefore, we have $\overline{v_1}=\beta_1(\overline{y},\overline{z})$ for some $\beta_1(y,z)\in R[y,z]$ with $\beta_1(y,z)\notin xB$. Then $v_1=\beta_1(y,z)+x^{t_2}v_2$ for some $v_2\notin xB$ and $t_2\geq 1$. From (\ref{e11}) we get, 
\begin{equation}\label{e13}
\psi(y,z)= x^{t_1}\beta_1(y,z) + x^{t_1+t_2} v_2.
\end{equation} 
Applying $D$ on \eqref{e13} we have 
\begin{multline}\label{e14}
x^m (\partial_z \psi(y,z))  + x^d f(x,z) (\partial_y \psi(y,z)) \ = x^{m+t_1} (\partial_z \beta_1(y,z)) + x^{d+t_1} f(x,z) (\partial_y \beta_1(y,z))\\+ x^{t_1+t_2} D(v_2).
\end{multline}  
Suppose $t_1+t_2 > d$, then using the facts that $\deg_z \partial_z \varphi(z) < \deg_z \varphi(z)$ and $\deg_y (\partial_y \psi(y,\overline{z}))$ $< \deg_y \psi(y,\overline{z})$ we get a similar contradiction as above. Thus $t_1+t_2 \leq d$. 
Proceeding such a way we obtain 
\begin{equation}\label{e15}
	\psi(y,z)=x^{t_1} \beta_1(y,z)+ x^{t_1+t_2} \beta_2(y,z)+\cdots+x^{t_1+\cdots+t_n}\beta_n(y,z)+x^{d} v
\end{equation}
where $t_1< t_1+t_2<\cdots<t_1+\cdots+t_n<d$ and $v\notin xB$.
Now applying $D$ on \eqref{e15}, we have 
\begin{equation}\label{e16}
Dv= f(x,z) (\partial_y \psi(y,z)) +xb_1
\end{equation} 
for some $b_1 \in B$. Since $f(x,z) (\partial_y \psi(y,z)) \notin xB$ we have $\overline{v} \notin \ker \overline{D}$.
 Hence $\overline{v}$ is algebraically independent over $R[\overline{y},\overline{z}]$.
 Now, note that $B \subseteq R[x,x^{-1}, y,z,v]$, $\overline{v}$ is algebraically independent over $R[\overline{y},\overline{z}]$ and \eqref{e19} and \eqref{e15} hold. Therefore, by \thref{lem1}, we have $B=R[x,y,z,v]$ where the generators satisfy the relations 
 \begin{align}
 	 &F(x,z)=x^{m-d}y\\
 	 &G(x,y,z)= x^dv
 \end{align}
such that \eqref{e8} holds, that means $F(0,Z)$ is irreducible in $R[z]$. Further, $G(x,y,z):=\psi(y,z)-x^{t_1} \beta_1(y,z)- x^{t_1+t_2} \beta_2(y,z)- \cdots - x^{t_1+\cdots+t_n}\beta_n(y,z)$ with $G(0,y,z)=\psi(y,z)$ which is almost monic in $y$ as a polynomial in $K[y,z]$. Therefore, by \thref{lmm3} the assertion follows.
 \end{proof}

\begin{rem}
\rm Note that, in Theorems \ref{DD} and \ref{char fac dd}, if $R=k$ then the condition that $k[\overline{z}]$ is factorial is redundant. 
\end{rem}

\end{document}